\documentclass[10pt,onecolumn,aps,prd,floats,floatfix,superscriptaddress,nofootinbib]{revtex4}
\usepackage{amsmath,amsfonts}
\usepackage[dvips]{graphicx}
\usepackage{color}
\usepackage{epsfig}
\begin{document}
\title{Variational Minimization on String-rearrangement Surfaces,\\ Illustrated by an Analysis of the Bilinear Interpolation.}
\author{Daud Ahmad}
\thanks{daudahmadpu@yahoo.com}
\affiliation {Department of Mathematics, University of the Punjab,Lahore, Pakistan}
\author{Bilal Masud}
\thanks{bilalmasud.chep@pu.edu.pk}
\affiliation {Center for High Energy Physics, University of the Punjab, Lahore, Pakistan}

\begin{abstract}
In this paper we present an algorithm to  reduce  the area of a surface spanned by a finite number of boundary curves  by initiating a variational improvement in  the surface. The ansatz we suggest consists of original surface plus a variational parameter $t$ multiplying the  numerator $H_{0}$ of mean curvature function defined over the surface. We point out that the integral of the square of the mean curvature with respect to the surface parameter becomes a polynomial in this variational parameter. Finding a  zero, if there is any, of this polynomial would end up at the same (minimal) surface as obtained by minimizing more complicated area functional itself. We have instead minimized this polynomial. Moreover, our minimization is restricted to a search in the class of all surfaces allowed by our ansatz.  All in all, we have not yet  obtained the exact minimal but we do reduce the area for the same fixed boundary. This reduction is significant for a surface (hemiellipsoid) for which we know the exact minimal surface. But for the bilinear interpolation spanned by four bounding straight lines, which can model the initial and final configurations of re-arranging strings, the decrease remains less than 0.8 percent of the original area. This may suggest that bilinear interpolation is already a near minimal surface.
\end{abstract}
\maketitle
\section{Introduction}
Variational methods are one of the active research areas of the optimization theory \cite{JXWD}.  A variational method tries to find the best values of the parameters in a trial function   that optimize, subject to some algebraic, integral or differential constraints, a quantity dependant on the ansatz. A simple example of such a problem may be to find the curve of shortest length connecting two points. The solution is a straight line between the points  in case of no constraints and simplest metric, otherwise possibly many solutions may exist depending on the nature of constraints. Such solutions are called geodesics \cite{PeyerCohn2008,improvmentSP,Li2013}.   One of the related problems is finding the path of stationary optical length connecting two points, as the Fermat's principle says that the rays of light traverse such a path.   Another related problem is a \textit{Plateau problem}  \cite{Osserman1986,Nitsche1989} which is   finding the surface with minimal area enclosed by a given curve.  This problem is named after the blind Belgian physicist Joseph Plateau, who demonstrated in 1849 that a minimal surface can be obtained by immersing a wire frame, representing the boundaries, into soapy water. The Plateau problem attracted mathematicians like Schwarz \cite{Schwarz} (who discovered D (diamond), P (primitive), H (hexagonal), T (tetragonal) and CLP (crossed layers of parallels) triply periodic surfaces), Riemann \cite{Osserman1986}, and Weierstrass \cite{Osserman1986}. Although mathematical solutions for specific boundaries had been obtained for  years, but it was not until 1931 that the American mathematician Jesse Douglas \cite{Douglas}   and the Hungarian Tibor Rad\'o \cite{Rado} independently proved the existence of a minimal solution for a given simple  closed curve. Their methods were quite different. Douglas \cite{Douglas} minimized a functional now named as Douglas-Dirichlet Integral. This is easier to manage but has the same extremals in an unrestricted search \cite{Monterde2004} as the area functional.  Douglas results  held for arbitrary simple closed curve, while Rad\'o \cite{Rado} minimized the energy.  The work of Rad\'o was built on the previous work of R. Garnier \cite{Garnier} and held only for rectifiable simple closed curves. Many results  were obtained in subsequent years, including  revolutionary achievements of L. Tonelli \cite{Tonelli},  R. Courant \cite{Courant1} \cite{Courant2}, C. B. Morrey \cite{Morrey1} \cite{Morrey2}, E. M. McShane \cite{Shane}, M. Shiffman \cite{Shiffman}, M. Morse \cite{MT}, T. Tompkins \cite{MT}, Osserman \cite{Osserman1970}, Gulliver \cite{Gulliver} and Karcher \cite{Karcher} and others.

In addition to finding (above mentioned) alternative functionals, the search can be limited to a certain class of surfaces.  A widely used such restriction is to search among all B\'{e}zier surfaces with the given boundary. B\'{e}zier models are widely used in computer aided geometric design (CAGD) because of their  suitable geometric properties.  For a  control net $\mathbf{P}_{ij}$ of a  two dimensional parametric B\'{e}zier surface is given by
\begin{equation}\label{beziersurface}
  \mathbf{x}(u,v)=\sum\limits_{i=0}^{n}{\sum\limits_{j=0}^{m}{B_{i}^{n}\left( u \right)B_{j}^{m}\left( v \right)}}\text{ }\mathbf{P}_{ij},
\end{equation}
where $u,v$ are the parameters, $B_{i}^{n}\left( u \right) = \binom {n}{i}{{u}^{i}}{{\left( 1-u \right)}^{n-i}}$,  the Bernstein polynomials of degree $n$ and $\binom {n} {i}= \frac{n!}{i!\left( n-! \right)}$,  binomial coefficients and $D = [0,1] \times [0,1]$. The minimal B\'{e}zier surfaces as an example of the extremal of  discrete version of Dirichlet functional may be found in the Monterde work \cite{Monterde2004},  a restricted Plateau-B\'{e}zier problem defined as the  surface of minimal area among all B\'{e}zier surfaces with the given boundary. A use of Dirichlet method and the extended bending energy method to obtain an approximate solution of Plateau-B\'{e}zier problem may be seen in work by Chen et al \cite{Chen2009}. Another restriction may be to find a surface in the parametric polynomial form as it can be seen in the ref.~\cite{Xu2010} that finds a class of quintic parametric polynomial minimal surfaces. B\'{e}zier surfaces exactly deal with the case that the prescribed borders are polynomial curves. A more general case of borders is taken in  ref.~\cite{Hao2012}  that study  the Plateau-quasi-B\'{e}zier problem which includes the  case when the boundary curves are catenaries and circular arcs.   The Plateau-quasi-B\'{e}zier problem is related to the quasi-B\'{e}zier surface with minimal area among all the quasi-B\'{e}zier surfaces with prescribed border. They minimize the Dirichlet functional in place of original area functional.

An emerging use of minimal surfaces in physics is that in string theories. A classical particle travels a geodesic with least distance whereas a classical string is an entity which traverses a minimal area. Amongst the string theories used in physics, two are worth mentioning. One is the theory of quantum chromodynamics (QCD) strings that model the gluonic field confining a quark and an antiquark within a meson. (The gluonic field connecting three quarks, within a proton or neutron, is modeled through Y-shaped strings. For a system composed of more than three quarks, minimization of the total length of a string network with only Y-shaped junctions may be a non-trivial  Steiner-Tree Problem \cite{JMR}).  In the other string theory (or theories) string vibrations are supposed to generate different elementary particles of the present high energy physics. Quite often string theories need a surface spanning the boundary composed of curves either connecting particles or describing the time evolution of particles. An important case can be a fixed boundary composed of four external curves. A common application of this boundary can be the time evolution of a string parameterized by $\sigma$ or $\beta$ \cite {SY} variable; the time evolution itself is parameterized by the symbol $\tau$, the proper time  of relativity.  In this case two bounding curves parameterized by the respective $\sigma$ or $\beta$ represent the initial and final configurations of a string, and the other two curves (parameterized by the respective $\tau$ variables) describe the time evolution of the two ends of a string.

String theories take action to be proportional to area. Combining this with the classical mechanics demand of the least action, minimal surfaces spanning the corresponding fixed boundaries get their importance. For example, see eq. 13 of ref. \cite{SY} for the Nambu-Goto ansatz for the minimal surface area and compare it with eqs. \eqref{BIEu} and \eqref{BIEv} below, along with ref.  \cite {AOV} for Nambu-Goto strings. Also relevant is the use in ref. \cite{BV} of Wilson minimal area law (MAL) to derive the quark antiquark potential in a certain approximation. A surface spanned by such a boundary is in space-time of relativity. An ordinary 3-dimensional spatial surface can span a boundary composed of two 3-dimensional curves connecting four particles and two other curves connecting the same four particles in a re-arranged (or exchanged) clustering; see for example Fig. 2 of ref. \cite{GP} and Fig. 5 of ref. \cite{GLW}. An explicit expression of such a spanning surface can be found in eqs. 3, 4 of ref. \cite{BF1} and eq. 22  of ref. \cite{FGM}. This is a bilinear interpolation in ordinary 3-dimensional space and is similar to the linear  interpolations in above mentioned eq. 13 of ref. \cite{SY}, eq. 4.7 of ref. \cite{BMP} and eq. 3.4 of ref. \cite{BV}.  Ref. \cite{BMP} clarifies that such a surface is used as a {\em replacement} to the exact minimal surfaces  for the corresponding boundaries; see section \ref{basics} below for a minimal surface in the differential geometry.  Even non-minimal surfaces have some usage in the mathematical modeling of quantum strings because 1) in contrast to classical strings, quantum strings can have any action and hence area as described by the path integral version of the quantum mechanics (see eq. 1 of   \cite{JM} ) and 2) any surface spanning a boundary composed of quark lines (or quark connecting lines) corresponds to a physically allowed (gauge invariant) configuration of the gluonic field between these quarks; compare the non-minimal surface of Fig. 10.5 of ref. \cite{MCz}  with the minimal surface for the same boundary in Fig. 10.1 of the same ref. \cite{MCz}.  But it cannot be denied that minimal surfaces are the most important of the spanning surfaces even in quantum theories. For example,  the relation in eq. 1.14 of ref. \cite{BCP}   between an area and an important quantity (termed {\em Wilson loop})  related to the potential between a quark and antiquark connected by a QCD (gluonic) string is valid only if the area is of the minimal surface. (Though above mentioned eq. 1 of ref. \cite{JM}   relates the Wilson loop to a ``sum over all surfaces of the topology of rectangle bounded by the loop" implying that each spanning surfaces has some contribution in the Wilson loop,  the minimal surface must contribute most.)  Thus it is worth pointing out that the non-minimal linearly or bilinearly interpolating surfaces  can replace minimal surfaces, can be effectively used as minimal surfaces or share some features in common with minimal surfaces;  text just before  eq. (1.15) of the above mentioned   ref.  \cite{BCP} relates them, up to  non-relativistic   1/(mass square) order, to the minimal surfaces. The purpose of the present paper is explore further this  ``effective usability" or "sharing common features with minimal surfaces" of linearly or bilinearly interpolating surfaces. Before starting a description of our work, we want to 1) state the common feature we have chosen. This is the fractional reduction possible in the area for a fixed boundary; for an exact minimal surface this quantity is zero (at least for a small neighbourhood).
For reducing area we use the variational area reduction, outlined in sect. \ref {tvi}, to our specific bilinear interpolation described in sec. \ref{ssfbc}. Moreover, we   2) point out that the bilinear interpolations used in string-theories-related works of physics are also used in the emerging discipline of the computer aided geometric design (CAGD) and hence the usefulness of the present paper extends to above mentioned CAGD along with   physics and the differential geometry;   as much as bilinear interpolations are near or related to minimal surfaces their study sheds some light on the above mentioned Plateau problem of the differential geometry itself.

Computer aided geometric design (CAGD) \cite{Farin}, \cite{Faux}, mentioned above, arose  when mathematical descriptions of shapes facilitated the use of computers to process data and analyze related information. In the 1960s, it became possible to use computer control for basic and detailed design enabling utilization of a mathematical model stored in a computer instead of the conventional design based on drawings. The term geometric modeling is  used to characterize the methods used in describing the geometry of an object.  Over the years, various schemes were developed with a view to achieve this abstraction.  S. A. Coons \cite{Farin}, \cite{Faux}  introduced the Coons patch in 1964. The Coons patch approach is based on the premise that a patch can be described in terms of four distinct boundary curves. Thus a Coons patch can be a worth analyzing surface spanning a fixed set of boundary curves. This is simple when the number of bounding curves is four. For a surface spanned by an arbitrary $N$-number of curves, it is still possible to find a Coons patch that is spanned by a boundary of four analytical curves by combining, as for example the way we did in ref.~\cite{dabm2013}, these $N$-number of curves into four groups and then joining these curves in each of four groups into a single analytic curve. This joining let us use eq.~\eqref{CPR} to write the Coons patch spanned by $N>4$-number of curves which may then be used to find the associated  minimal surface by the ansatz eq. \eqref{VS1}. Using that formalism our technique can be applied to any number of curves, we have implemented it in full though numerical implementation has been limited to five straight lines. Ref.\cite{Monterde2006} points out that Coons patch can be considered a special case of the above mentioned B\'{e}zier surface. For us, Coons patch (see eq. \eqref{CPR} below) is relevant because the above mentioned   bilinear interpolations (see also eq. \eqref{BIE} below) we basically study in this paper are a special case of Coons patch \cite{Farin}. Coons patch analysis is an active area of research and has seen enormous development during recent years. But most, if not all, of the work on it has been limited to its geometric descriptions and visualization and to interactive mathematical experiments with it; it has not been analyzed from the view of differential geometry and that is also what we aim to do in this paper though we actually study only its special case of a bilinear interpolation. In trying to judge how close it is to being a minimal surface, we see how much its area can be reduced through our variational minimization. To carry out his minimization, we also had to restrict our surface search to surfaces of the form of eq.~\eqref{VS1} below. This restriction can be compared to the more common above mentioned restriction to the B\'{e}zier surfaces. As for minimization, we have used the mean square mean curvature of our eq.~\eqref{rms}.

The paper is organized as follows. In  sections \ref{basics} and \ref{ssfbc} we present basic  definitions and  constructions related to surfaces spanned by fixed boundary curves. In the next section  \ref{tvi} we  present an algorithm to  reduce  the area of a surface spanned by a finite number of boundary curves by introducing a variational improvement in  a surface. Then in section \ref{hessfal} we apply this technique to reduce  the area of a non-minimal surface spanning a boundary for which the minimal surface is known - namely hemiellipsoid  eq. \eqref{hemiellipsoid}, to make sure the efficiency of the algorithm given by eq. \eqref{VS1} and   above mentioned bilinear interpolation spanned by four bounding lines  for which the corresponding minimal surface is not known. Based on this comparison, we comment on the possible status of bilinear interpolation as an approximate minimal surface. The last section \ref{conclusion} presents results, final remarks and mentions possible future developments.
\section{Differential Geometry of Minimal Surfaces}  \label{basics}
 In the optimization problem we aim for here, we eventually try to find a surface of a known boundary that has a least value of area. Area is evaluated by the area functional:
\begin{equation} \label{af}
A(\mathbf{x})=\int\int_{D}\left| \mathbf x_{u}(u,v)\times \mathbf x_{v}(u,v) \right|du dv,
\end{equation}
where $ D \subset R^{2} $ is a domain over which the surface $ \mathbf x(u,v) $ is defined as a map, with the boundary condition $\mathbf x(\partial D)=\Gamma$ for $ 0\leq u\leq1 $ and $ 0\leq v\leq1 $,   $\mathbf x_{u}(u,v)$ and $\mathbf x_{v}(u,v)$ being partial derivatives of $ \mathbf x(u,v) $ with respect to $u$ and $v$ .
It is known  \cite{docarmo} that the first variation of $A(\mathbf x)$ vanishes everywhere if and only if the mean curvature $H$ of $\mathbf x(u,v)$ is zero everywhere in it. Thus a surface of least area is also a surface of least (zero) $rms$ mean curvature spanning the given boundary. This means we can aim for the same surface using the condition of the least means square mean curvature in place of the condition of the least area. This is helpful as, unlike area, the ms mean curvature has not a square root in its integrand. For a locally parameterized surface $ \mathbf x= \mathbf x(u,v) $, the mean curvature $H$  may be given by
\begin{equation}
H =\frac{G e -2 F f+ E g}{E G -F^{2}}, \label{mc}
\end{equation}
where
\begin{equation} \label{ffm}
E =\left\langle \mathbf x_{u},\mathbf x_{u}\right\rangle,\hspace{0.5cm}F =\left\langle \mathbf x_{u},\mathbf x_{v}\right\rangle,\hspace{0.5cm} G =\left\langle \mathbf x_{v},\mathbf x_{v}\right\rangle,
\end{equation}
are the first fundamental coefficients and
\begin{equation} \label{sfm}
e =\left\langle \mathbf N,\mathbf x_{uu}\right\rangle, \hspace{0.5cm} f =\left\langle \mathbf N,\mathbf x_{uv}\right\rangle, \hspace{0.5cm} g =\left\langle \mathbf N,\mathbf x_{v}\right\rangle,
\end{equation}
are the second fundamental coefficients with
\begin{equation} \label{un}
\mathbf N(u,v) =\frac{\mathbf x_{u} \times \mathbf x_{v}}{\left| \mathbf x_{u} \times \mathbf x_{v}\right|},
\end{equation}
being the unit normal to the surface $\mathbf x(u,v) $.
The root mean square root  of the mean curvature $H(u,v)$, for $ 0\leq u\leq1
$ and $ 0\leq v\leq1 $ denoted by $\mu$ is given by the following
expression,
\begin{equation} \label{rms}
\mu = \left({\int^{1}_{0} \int^{1}_{0}H^{2}\hspace{0.2cm}dudv}\right)^{1/2}.
\end{equation}
For a minimal surface \cite{docarmo}, \cite{Goetz} the mean curvature \eqref{mc} is identically zero. For minimization we use only the numerator part of mean curvature $H$ given by \eqref{mc}, as done in ref. \cite{BCDH} following ref. \cite{Osserman1986} who writes that   ``for a locally parameterized surface, the mean curvature vanishes when the numerator part  of the mean curvature is equal to zero".  We call this numerator part $H_{0}$  as the  $rms$ curvature of the initial surface   $\mathbf x_{0} (u,v)$   to be used in the ansatz  eq. \eqref{VS1} to get first order variationally improved surface $\mathbf x_{1}(u,v)$ of lesser area.   This process could be continued as an iterative process until a minimal surface is achieved.   But due to complexity of the calculations required for obtaining the second order improvement $\mathbf x_{2}(u,v)$, we have been able to calculate  the first order $\mathbf x_{1}(u,v)$ only. The numerator part $H_0$ is denoted by
 \begin{equation}
 H_{0}=e_{0} G _{0}-2F_{0} f_{0}+g_{0} E_{0} \label{NMC},
\end{equation}
where $E_0, F_0, G_0, e_0, f_0 $ and $g_0$ denote the fundamental magnitudes given by eqs. \eqref{ffm} and \eqref{sfm}, with $N_0(u,v)$ being the unit normal given by eq.  \eqref{un} to the initial surface   $\mathbf x_{0} (u,v)$. We call the root mean square $(rms)$ of this $H_0$, for $ 0\leq u\leq1 $ and $ 0\leq v\leq1 $, as $\mu_{0}$. That is,
\begin{equation}
\mu_{0}= \left({\int^{1}_{0} \int^{1}_{0}H_{0}^{2}\hspace{0.2cm}dudv}\right)^{1/2}. \label{rms0}
\end{equation}
In the notation of eqs \eqref{mc} to \eqref{sfm} eq. \eqref{af} becomes, for $\mathbf x_0 (u,v)$,
\begin{equation}
A_{0}= \int^{1}_{0} \int^{1}_{0} \sqrt{E_{0} G_{0}-F_{0}^{2}} \hspace{0.2cm} dudv. \label{area0}
\end{equation}

\section{Bilinear Starting Surface Spanned by Fixed Boundary Curves} \label{ssfbc}
 For a minimal (or, more precisely, a stationary) surface, we have to solve the differential equation obtained by setting the mean curvature $H$ given by eq. (\ref{mc}) equal to zero for each value of the two parameters, say, $u$ and $v$ parameterizing a surface spanning the fixed boundary. In this section, our purpose is to describe a starting surface bounded by the skew quadrilateral which is composed of four arbitrary straight lines connecting four corners $\mathbf x(0,0), \mathbf x(0,1), \mathbf x(1,0)$ and $\mathbf x(1,1)$; in the next section we report the variational improvement to this start aimed towards minimizing the surface evolving from the starting surface of this section.  Ref. \cite{BF2} also includes a preliminary effort to variationally improve the surface bounded by four straight lines towards being a minimal surface. The algorithm used for this variational improvement applies to a wider class of surfaces. Accordingly, now we point out a class of surfaces, namely Coons patch, that includes surfaces bounded by four straight lines: Let $\mathbf c_{1}(u), \mathbf c_{2}(u), \mathbf d_{1}(v)$ and $\mathbf d_{2}(v)$  be four given arbitrary curves defined over the parameters $u,v \in \left[0,1\right]$ . For  $\mathbf c_{1}(u)=\mathbf x(u,0)$, $\mathbf c_{2}(u)=\mathbf x(u,1)$, $\mathbf d_{1}(v)=\mathbf x(0,v)$ and $\mathbf d_{2}(v)=\mathbf x(1,v)$, blending functions $f_{1}(u)$, $f_{2}(u)$, $g_{1}(v)$ and $g_{2}(v)$ satisfying the conditions that $f_{1}(u)+f_{2}(u)=1$, $g_{1}(v)+g_{2}(v)=1$ for non-barycentric combination of points and $f_{1}(0)=g_{1}(0)=1, f_{1}(1)=g_{1}(1)=0$ in order to actually interpolate $\mathbf x(0,0), \mathbf x(0,1), \mathbf x(1,0)$ and $\mathbf x(1,1)$,  the following equation defines Coons patch:
\begin{equation}
\begin{split}
\mathbf x(u,v)& =
\left[\begin{array}{cc} f_{1}(u) & f_{2}(u)  \end{array} \right]
\left[\begin{array}{cc} \mathbf x(0,v)\\ \mathbf x(1,v) \end{array} \right] +
\left[\begin{array}{cc} \mathbf x(u,0) & \mathbf x(u,1) \end{array} \right]
\left[\begin{array}{cc} g_{1}(v)\\g_{2}(v) \end{array} \right]  \\ & \qquad
- \left[\begin{array}{cc} f_{1}(u) & f_{2}(u) \end{array} \right]
\left[\begin{array}{cc} \mathbf x(0,0) & \mathbf x(0,1) \\ \mathbf x(1,0) & \mathbf x(1,1)  \end{array} \right]
\left[\begin{array}{cc} g_{1}(v)\\g_{2}(v) \end{array} \right].
\end{split} \label{CPR}
\end{equation}

As a special case of the above, we consider a Coons patch for which all the three terms in eq.  \eqref{CPR} are equal,  so that this equation reduces to the following form:
\begin{equation}
\begin{split}
\mathbf x(u,v)& =
\left[\begin{array}{cc} f_{1}(u) & f_{2}(u) \end{array} \right]
\left[\begin{array}{cc} \mathbf x(0,0) & \mathbf x(0,1) \\ \mathbf x(1,0) & \mathbf x(1,1)  \end{array} \right]
\left[\begin{array}{cc} g_{1}(v)\\g_{2}(v) \end{array} \right].
\end{split} \label{BIE}
\end{equation}

 The boundary spanned by lines connecting the points $\mathbf x(0,0)$, $\mathbf x(0,1)$, $\mathbf x(1,0)$ and $\mathbf x(1,1)$  with linear blending functions
  \begin{equation}\label{lbf}
     f_{1}=1-u, \hspace{.5cm}  f_{2}=u , \hspace{.5cm}  g_{1}=1-v, \hspace{.5cm}   g_{2}=v,
  \end{equation}
  in eq. \eqref{BIE} can represent a  time evolution of a string or, alternatively, a re-arrangement of a one set of two strings to the only possible other re-arranged set (of two strings) connecting the same two particles and two antiparticles.  (It is to be noted that a string connects only a particle with antiparticle. This constraint allows only two string arrangements for a system composed of two particles and two antiparticles.) Above $\mathbf x(u,v)$ eq. \eqref{BIE}   spanning this boundary is a surface that is needed in many models of string re-arrangements from one of these configurations to the other with the  particle positions $\mathbf x(0,0) \equiv 1$ and $\mathbf x(1,1) \equiv 2$ and anti- particle positions $\mathbf x(1,0) \equiv \bar{3}$ and $\mathbf x(0,1) \equiv \bar{4}$. In this paper we reduce the area of a quadrilateral, using above  linear blending functions  and particle positions. This gives
\begin{equation}
\mathbf x_{u}(u,v) =\left(1-v\right) \mathbf r_{1 \bar{3}} -v \mathbf r_{2 \bar{4}}, \label{BIEu}
\end{equation}
and
\begin{equation}
\mathbf x_{v}(u,v) =\left(1-u\right) \mathbf r_{1 \bar{4}} -u \mathbf r_{2 \bar{3}}, \label{BIEv}
\end{equation}
 as partial derivatives $w. r. t.$ $ u$ and $v$
 with the following corners:
\begin{equation}
	\mathbf x(0,0)=\mathbf r_{1}, \hspace{.5cm}x(1,1)=\mathbf r_{2}, \hspace{.5cm}\mathbf x(1,0)=\mathbf r_{\bar{3}}, \hspace{.5cm} \mathbf x(0,1)=\mathbf r_{\bar{4}},
\end{equation}
($\mathbf x(u,v)$ is our starting surface spanning four straight lines.) For real scalars $r$ and $d$, we consider two types of configurations of the four corners: ${\text {ruled}}_1$ and ${\text {ruled}}_2$. For ${\text {ruled}}_1$ we choose
\begin{equation} \label{ruled1}
	\mathbf r_{1}=(0,0,0),  \hspace{.5cm} \mathbf r_{2}=(r,d,0),  \hspace{.5cm} \mathbf  r_{\bar{3}}=(0,d,d),  \hspace{.5cm} \mathbf r_{\bar{4}}=(r,0,d).
\end{equation}
The mapping from $(u,v)$ to $(x,y,z)$ in this case is
\begin{equation}
x(u,v)=r\left(u+v-2u v \right), \hspace{.5cm} y(u,v)=v d, \hspace{.5cm} z(u,v)=u d.
\end{equation}
For ${\text {ruled}}_2$ we choose
\begin{equation} \label{ruled2}
	\mathbf x(0,0)=\mathbf r_{1}, \hspace{.5cm} \mathbf x(1,1)=\mathbf r_{2},  \hspace{.5cm} \mathbf x(0,1)=\mathbf r_{\bar{3}}, \hspace{.5cm} \mathbf x(1,0)=\mathbf r_{\bar{4}} .
\end{equation}
The mapping from $(u,v)$ to $(x,y,z)$ in this case is
\begin{equation}
x(u,v)=u r, \hspace{.5cm} y(u,v)=v d, \hspace{.5cm} z(u,v)=u d+v d(1-2u).
\end{equation}

\begin{figure}[htb!]
\begin{center}
\includegraphics[width=60mm]{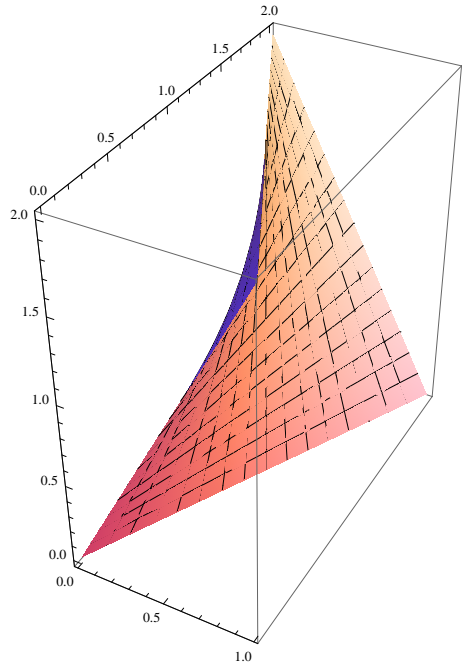}
\end{center}
\caption{The $ruled_{1}$  surface $(r=1, d=2) $with$x,y$ as the  horizontal plane and height along $z- axis$.} \label{figruled1}
 \end{figure}

\begin{figure}[htb!]
\begin{center}
\includegraphics[width=60mm]{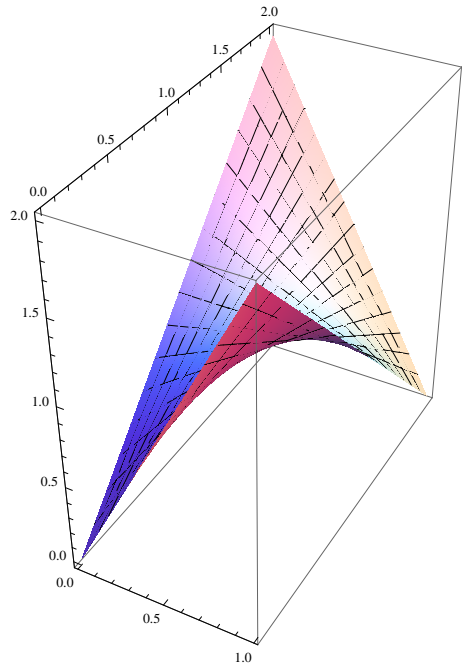}
\end{center}
\caption{The $ruled_{2}$ surface $(r=1, d=2)$ with $x,y$ as the  horizontal plane and height along $z- axis$.} \label{figruled2}
\end{figure}

These definitions are such that for $r=d$ the four position vectors $ \mathbf r_{1}, \mathbf r_{2},  \mathbf r_{3}$ and $ \mathbf r_{4}$ lie at the corners of a regular tetrahedron. Fig. \ref{figruled1} and Fig. \ref{figruled2} below are 3D graphs of surfaces called the hyperbolic paraboloids
  for a choice of corners given by  \eqref{ruled1} and \eqref{ruled2}.
The expression for the mean curvature, calculated using eq. \eqref{mc}, of our bilinear interpolations  is
\begin{equation}
-\frac{4 r^3 (2 u-1) (2 v-1)}{d \left(d^2+2 r^2 (2 (u-1) u+2 (v-1)  v+1)\right)^{3/2}},
\end{equation}
for the ${\text {ruled}}_1$ and
\begin{equation}
\frac{4 d r (2 u-1) (2 v-1)}{\left(d^2 (1-2 v)^2+2 r^2 (2 (u-1)  u+1)\right)^{3/2}},
\end{equation}
for the ${\text {ruled}}_2$.

The mean curvature for our starting surface is zero only for the $u=\frac{1}{2}$ line and the $v=\frac{1}{2}$ line, whereas for a minimal surface this should be zero for all values of $u$ and $v$. Below we describe our effort to improve our surface towards being minimal.
\section{A Technique For Variational Improvement} \label{tvi}

The area functional given by eq.~\eqref{af} is highly non-linear and  is difficult to minimize due to its high non-linearity. Douglas replaced it by the extremal-sharing Dirichlet functional
 \begin{equation}
   D \left(\mathbf{x} \right)=\frac{1}{2}\iint\limits_{D}{\left( \mathbf{x}_{u}^{2}+\mathbf{x}_{v}^{2} \right)}dudv=\frac{1}{2}\iint\limits_{D}{\left( E+G \right)}dudv,
\end{equation}
 to give his famous solution to Plateau problem. A list of other possibilities of such functionals can be found in ref.\cite{MUComparative,Chen2009}. The Dirichlet Integral is related to the area functional eq.~\eqref{af} by the following relation
\begin{equation}
  {{\left( EG-{{F}^{2}} \right)}^{\frac{1}{2}}}\le {{\left( EG \right)}^{\frac{1}{2}}}\le \frac{E\text{ }+G}{2}.
\end{equation}
Thus, for a  surface $\mathbf{x}(u,v)$, $A(\mathbf{x}) \leq D\left( \mathbf{x} \right)$. The equality of the two integrals holds only for an isothermal patch i.e. for which $E=G$ and $F=0$. Both the functionals are defined as the integrals of positive functions, thus they are bounded below and both the functions have a minimum for a compact domain. Thus finding a surface of minimal area is equivalent to solving variational problem of finding a surface with appropriate boundary conditions for which the integrals are minimum. Douglas suggested minimizing the Dirichlet integral that has the same extremal as the area functional. We suggest another functional that has the same extremal as the area integral. This is based on observation that for the extremal (minimal surface) of the area functional, the mean curvature is zero and hence an integral of the square of mean curvature would be least for this area. (This is because this integral eq.~\eqref{rms} is non-negative by construction and hence zero is its least value.) Thus a minimal surface is also an extremal of the $\mu^{2}$ (eq.~\eqref{rms}) along with being an extremal of the area functional. Now, as with Dirichlet integral, $\mu^{2}$ has no square root unlike the area integral. Others \cite{Monterde2004,MonterdeUgail2004} have converted Dirichlet integrals to a system of linear equations for inner control points in terms of known boundary control points. We can convert our $\mu_{1}^{2}$ eq.~\eqref{rms1} to polynomial in a variational parameter $t$ introduced through the ansatz eq.~\eqref{VS1}. In Monterde work \cite{Monterde2004}, a surface may be spanned by given control points as is the case with B\'{e}zier surface \cite{FD}. We are considering a surface that is spanned by the fixed boundary, though our straight line boundaries are in turn dictated by corner points. The Coons patch we are basing on is, according to ref.~\cite{Monterde2006}, is a special case of the B\'{e}zier surface eq.~\eqref{beziersurface}.

In contrast to the work mentioned in above references, we choose the minimization of $ms$  mean curvature  to reduce the area of a non-minimal surface $\mathbf x(u,v)$ in order to get a smooth variationally improved surface instead of Dirichlet integral. This $ms$ mean curvature functional is positive as the integrand is positive and is zero only for a minimal surface. Thus we try to find that value of variational parameter that makes this $ms$ mean curvature zero or the neighbouring value for which the resulting variational surface is minimal or has reduced area. These surfaces are spanned by a fixed boundary curve, as is the case with the hemiellipsoid eq.\eqref{hemiellipsoid} or the surfaces (eqs.~\eqref{CPR ruled1}) spanned by four boundary curves. The area reduction in the surface  bounded by a skew quadrilateral composed by four straight lines (see below eq.~\eqref{CPR ruled1}) is included in the section (see section \ref{hessfal}), whereas for a surface spanned by $N>4$-number of curves, we have developed (see ref.~\cite[eq.~16]{dabm2013}) a formalism  that groups these curves into four and then in each group these curves are joined using step-function representation (ref.~\cite[eqs.~24-26]{dabm2013}) into an analytic curve. Using that formalism we are able to write Coons patch out of it which can be used to find a variationally improved surface of reduced area by the ansatz eq.~\eqref{VS1}. The reduction scheme follows in the remaining  part of the present section.

We want to reduce area of a non-minimal surface $\mathbf x(u,v)$ using the expectation that reduced value of  $ms$  mean curvature, denoted by  $\mu^{2}$, in turn reduces the area $A$ of the surface $\mathbf x(u,v)$. As mentioned above, the $ms$ mean curvature $\mu^{2}(t)$ reduces to a polynomial in the variational parameter $t$ and can be solved for its minimum value as discussed above. Our scheme is to reduce the area of a surface $\mathbf x(u,v)$   given by eq.  \eqref{BIE}- a special case of  eq. \eqref{CPR},  spanned by a fixed boundary, by obtaining a variationally selected surface $\mathbf x_{1}(u,v)$ of lesser area. For the variational improvement in surface  \eqref{CPR},  we suggest an ansatz  essentially consisting of the original surface  $\mathbf x_{0}(u,v)$ of eq.~ \eqref{BIE}   plus a variational parameter multiplying the numerator of its mean curvature. In our notation it becomes
\begin{equation}\label{VS1}
\mathbf x_{1}(u,v,t)=\mathbf x_{0}(u,v)+t\,m(u,v) \mathbf k,
\end{equation}
where $t$ is our variational parameter and
\begin{equation}
m(u,v)=uv(1-u)(1-v)H_{0}, \label{vp}
\end{equation}
is chosen so that the variation at the boundary curves $u=0, u=1, v=0$ and $v=1$ is zero.  $\mathbf k$  is  a unit vector chosen such that it makes a small angle  with the normal to the original surface and  $H_{0}$, given by  \eqref{NMC}, is numerator of the initial mean curvature of the starting surface $\mathbf x_{0}(u,v)$. Calling the  fundamental magnitudes for $\mathbf x_{1}(u,v)$ as $E_{1}\left(u,v,t\right), F_{1}\left(u,v,t\right), G_{1}\left(u,v,t\right), e_{1}\left(u,v,t\right), f_{1}\left(u,v,t\right)$ and $g_{1}\left(u,v,t\right)$, the  area $A_{1}$ of the surface $\mathbf x_{1}\left(u,v,t\right)$  for $ 0\leq u\leq1 $ and $ 0\leq v\leq1 $  is given by
\begin{equation}\label{area1}
A_{1}= \int^{1}_{0} \int^{1}_{0} \sqrt{E_{1}G_{1}-F_{1}^{2}} \hspace{0.2cm} dudv.
\end{equation}
We denote the  numerator of mean curvature for $\mathbf x_{1}(u,v)$ eq. \eqref{VS1} as $H_{1}(u,v,t)$. It would have the following familiar expression
\begin{equation}
H_{1}(u,v,t) = E_{1} g_{1}-2F_{1} f_{1}+G_{1} e_{1}. \label{MC1}
\end{equation}
As $H_{1} ^{2} (u,v,t)$ is a polynomial in  $t$, with real coefficients $h_{i} (u,v)$, we  rewrite eq. \eqref{MC1} in the form
\begin{equation}
H^{2}_{1}(u,v,t)=\sum^{n}_{i=0} (h_{i} (u,v)) \hspace{0.1cm} t^{i}. \label{H_{1}}
\end{equation}

Here $n$ turns out to be  $10$; there being no higher powers of $t$ in the polynomials as it can be seen from the expression for $E_{1}(u,v,t)$, $F_{1}(u,v,t)$ and $G_{1}(u,v,t)$ which are quadratic in $t$ and $e_{1}(u,v,t)$, $f_{1}(u,v,t)$ and $g_{1}(u,v,t)$ which are cubic in $t$. Integrating (numerically if needed) these coefficients $w. r. t.$ $u$ and $v$ in the range $0\leq u,v\leq 1$ we get the following integral  for the mean square mean curvature
\begin{equation}\label{rms1}
\mu_{1}^{2} (t)= \int^{1}_{0} \int^{1}_{0} H_{1} ^{2} (u,v,t) \hspace{0.2cm} dudv  = t^{i} \int^{1}_{0} \int^{1}_{0} \sum^{n}_{i=0} (h_{i} (u,v)) \hspace{0.1cm} du dv.
\end{equation}
The expression in the parentheses on right hand side of above equation turns out to be a polynomial in $t$ of degree $n$, which can be minimized $w. r. t.$ $ t $ to find $t_{min}$. The resulting value of $t$ completely specify $new$ surface $\mathbf x_{1}(u,v)$. New surface $\mathbf x_{1}(u,v)$ is expected to have lesser area than that of original surface $\mathbf x_{0}(u,v)$.

 In order to see a geometrically meaningful (relative) change in area we calculate the dimension less area by dividing the difference of the (original) area of the Coons patch and the variationally decreased area by the original area.
\section{The Technique Applied to Hemiellipsoid and a surface spanned by Four Arbitrary Lines} \label{hessfal}

In this section we apply the technique introduced in the above section  \ref{tvi}  to reduce the area of two types of surfaces. In first instance we apply this technique to reduce the area of a non- minimal surface spanning a boundary for which the minimal surface is known namely hemiellipsoid  eq. \eqref{hemiellipsoid} whose boundary is an ellipse lying in a plane and thus minimal area in this case is that of the elliptic  disc. The reduction in area in this case makes sure the efficiency of the algorithm given by eq.  \eqref{VS1}. In the second example  we apply this technique to reduce the area of a bilinearly interpolating surface spanned by four boundary lines lying in different planes for which the corresponding minimal surface is not known.
\subsection{Hemiellipsoid-A Surface with Corresponding Known Minimal Surface}

\begin{figure}[htb!]
\begin{center}
\includegraphics[width=60mm]{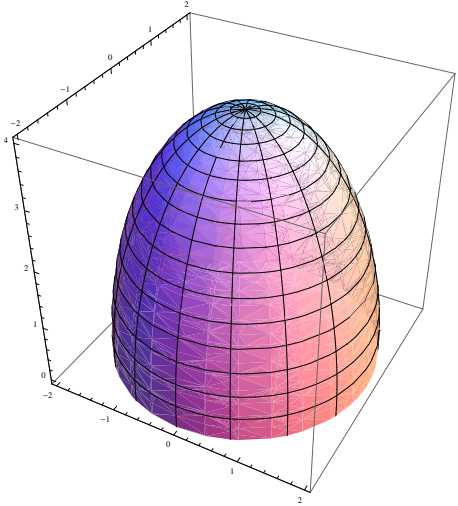}
\end{center}
\caption{A typical hemiellipsoid, initial non-minimal surface of which boundary is an ellipse in the $xy-plane$. } \label{hemiellipsoid1}
\end{figure}

We apply the technique introduced in the section  \ref{tvi} to the following surface $\mathbf x(u,v)$ namely hemiellipsoid  given by  eq. \eqref{hemiellipsoid} below along with linear blending functions eq. \eqref{lbf}, whose boundary is an ellipse. Simpler alternative of the above mentioned unit normal $\mathbf N(u,v)$, making a small angle  with it, in case of hemiellipsoid  eq. \eqref{hemiellipsoid} is found to be
\begin{equation}
\mathbf k=(0,0,1).
\end{equation}
A hemiellipsoid
\begin{equation}\label{hemiellipsoid}
    \mathbf x_{0} \left(u,v\right) = ( \sin u \cos v,b \sin u \sin v,c \cos u ).
\end{equation}
with $b$ and $c$ being constants and  $0\leq u\leq \frac{\pi}{2}$ and $0\leq v\leq 2 \pi$,  is a non-minimal surface  with its bounding curve  an ellipse in the $xy$-plane; see Fig. \ref{hemiellipsoid1}.
In this case we shall treat hemiellipsoid  as the initial non-minimal surface and the elliptical disc as the minimal surface for the given boundary, namely the ellipse. Thus, eq. \eqref{NMC} along with eqs.  \eqref{ffm} and \eqref{sfm} gives
\begin{equation}\label{mc0h}
   H_0=  -b c \sin ^3 u \left(\sin ^2 u \left(b^2 \cos ^2 v+c^2+\sin ^2 v\right)+\left(b^2+1\right) \cos ^2 u\right).
\end{equation}
The mean square mean curvature of beginning curvature given by eq. \eqref{rms0} takes the form
\begin{equation}
   \mu_{0}^{2}= {\int^{1}_{0} \int^{1}_{0}  b^2 c^2 \sin ^6 u \left(\sin ^2 u \left(b^2 \cos ^2 v+c^2+\sin ^2 v\right)+\left(b^2+1\right) \cos ^2 u\right)^2 \hspace{0.2cm}dudv}.
\end{equation}
The beginning or initial area of the Coons patch given by eq. \eqref{area0} takes the form in this case
\begin{equation}
   A_{0}=\int _0^1\int _0^1\sqrt{ \sin ^2 u \left(c^2 \sin ^2 u \left(b^2 \cos ^2 v+\sin ^2 v\right)+b^2 \cos  ^2 u\right)} \hspace{0.2cm}du dv.
\end{equation}
For Hemiellipsoid, $m(u,v)= (\pi/2-u)H_0 \, (0 \leq u\leq \pi/2)$ is the function that is zero at the boundary of the hemiellipsoid given by $u=\pi/2$. For $H_{0}$ from eq. \eqref{mc0h}  gives us expression for $m(u,v)$, thus in this case eq. \eqref{VS1} becomes
\begin{align}\label{vs1h}
\begin{split}
	\mathbf x_{1} (u,v,t) &=  (\sin  u \cos  v,b \sin  u \sin  v,c \cos  u-\frac{1}{16} b c t (\pi -2 u) \sin ^3 u (2 (b^2-2 c^2+1) \cos (2 u)-b^2  \cos (2 (u+v))-  (b^2-1) \\&\qquad \cos (2 (u-v))+2 b^2 \cos (2 v)+6 b^2+4 c^2+\cos (2 (u+v))-2 \cos (2 v)+6)).
\end{split}
\end{align}

\begin{figure}[htb!]
\begin{center}
\includegraphics[width=60mm]{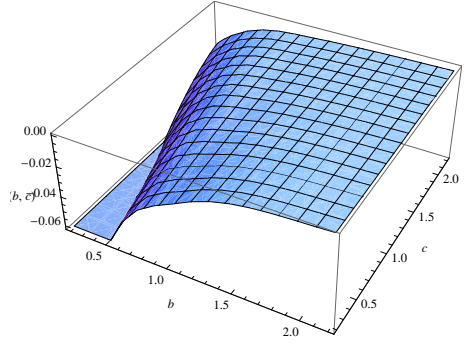}
\end{center}
\caption{Variation in  the parameter $t(b,c)$ as the semi-major and semi-minor axes $b$ and $c$ of the ellipse bounding the hemiellipsoid vary.} \label{hemiellipsoid2}
\end{figure}

Finding the   fundamental magnitudes $E_{1}\left(u,v,t\right), F_{1}\left(u,v,t\right), G_{1}\left(u,v,t\right), e_{1}\left(u,v,t\right), f_{1}\left(u,v,t\right)$ and $g_{1}\left(u,v,t\right)$  for the above surface  $\mathbf x_{1}(u,v)$  eq. \eqref{vs1h}, we can obtain the  area $A_{1}$  using eq. $\eqref{area1}$ and $H_{1}(u,v,t)$  using eq. \eqref{MC1} and after performing the integrations mentioned in eq. \eqref{rms1}, the  mean square curvature $\mu_{1}^{2}(t)$ for $\mathbf x_{1}(u,v)$ can be calculated. These are the similar details as given below for the non-minimal surface spanned by $4-$ non-coplanar lines. They have not been included for this ``first instance" but rather included for the ``second example"  because the formalism is well illustrated by this ``second example".  Also, that these expressions are too lengthy to be presented.  For chosen  values of $b$ and $c$ we can  generate a table of their values  within the  range $0\leq u\leq \frac{\pi}{2}$ and $0\leq v\leq 2 \pi$. For our purpose we took  $0\leq b,c\leq2$ with a step size $0.2$  and $0\leq u\leq \frac{\pi}{2}$ and $0\leq v\leq 2 \pi$ , yielding a table of values.  Interpolation surface for the corresponding minimum values $t(b,c)$ as a function of $b$ and $c$ is given by Fig.\ref{hemiellipsoid2}.

\begin{figure}[htb!]
\begin{center}
\includegraphics[width=60mm]{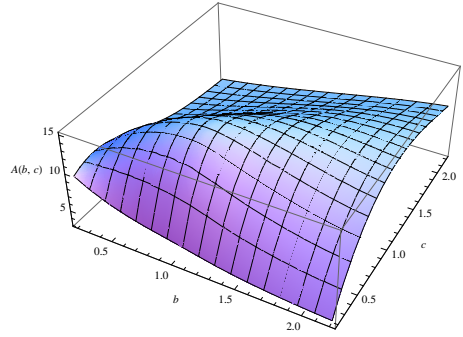}
\end{center}
\caption{The dimensionless decrease in area $A(b,c)$ of  hemiellipsoid
as a function of semi-major and semi-minor axes $b$ and $c$ of the ellipse bounding the hemiellipsoid.} \label{hemiellipsoid3}
\end{figure}

In this case the dimensionless decrease  $p$ in area for different values of $b$ and $c$ is  $0\leq p \leq 15$ that may be seen from the Fig. \ref{hemiellipsoid3}.

\subsection{Surface spanned by four arbitrary boundary lines}

\begin{figure}[htb!]
\begin{center}
\includegraphics[width=60mm]{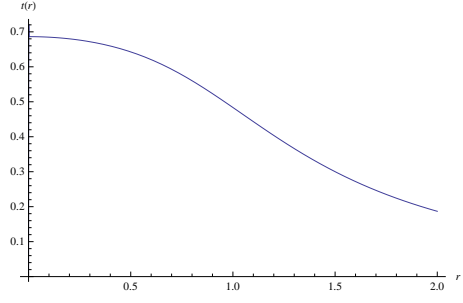}
\end{center}
\caption{The variation in parameter $t(r)$ depends on the variation of real scalar  $r$, for the  skew quadrilaterals ${\text {ruled}}_1$  bounded by  four arbitrary straight lines connecting four corners $ x(0,0), x(0,1),  x(1,0)$ and $ x(1,1)$.
} \label{interpolation1}
\end{figure}

Now we apply the technique introduced in the section  \ref{tvi} to the eq. \eqref{BIE} along with linear blending functions  eq. \eqref{lbf} for a surface $\mathbf x(u,v)$  whose boundary is composed of  four  straight lines connecting four arbitrary corner points $\mathbf x(0,0), \mathbf x(0,1), \mathbf x(1,0)$ and $\mathbf x(1,1)$. For its corners we choose  the configuration eq. \eqref{ruled1}, for a selection of integer values of r and d. The results  for the configuration \eqref{ruled2} have not been included as they are similar  to those for the configuration \eqref{ruled1}.  We found that the above mentioned simpler alternative of the unit normal $\mathbf N(u,v)$, making a small angle  with it, in case of  configuration eq. \eqref{ruled1} is
\begin{equation}
\mathbf k=(-1,0,0).
\end{equation}
Inserting values of blending functions and boundary points in the eq. \eqref{BIE} we find
\begin{equation}
	\mathbf x_{0}(u,v) =\left(r (u+v-2 u v),v d,u d \right), \label{CPR ruled1}
\end{equation}
with fundamental magnitudes having the expressions as
\begin{equation}
	E_{0}= d^{2}+r^{2} (1-2v)^{2}, \hspace{0.25cm} F_{0} = r^{2} (1-2u) (1-2v), \hspace{0.25cm}  G_{0} =d^{2}+r^{2} (1-2 u)^{2},
\end{equation}
\begin{equation}
	e_{0}= 0,   \hspace{0.25cm}  f_{0} = 2 d^{2} r,  \text{and} \hspace{0.25cm}   g_{0} = 0.
\end{equation}
Thus, eq. \eqref{NMC} gives
\begin{equation}
	H_{0}=-4 d^{2} r^{3} (-1+2 u) (-1+2 v). \label{mc0}
\end{equation}
 The root mean square ($rms$) of beginning curvature given by eq. \eqref{rms0} takes the form
\begin{equation}
	\mu_{0}=\frac{4 d^{2} r^{3}}{3}.
\end{equation}
The beginning or initial area of the Coons patch given by eq. \eqref{area0} takes the form in this case
\begin{equation}
  A_{0}=\int^{1}_{0} \int^{1}_{0}d \sqrt{d^2+2 r^2 \left(2 u^2-2 u+2 v^2-2 v+1\right)} \hspace{0.2cm}du dv. \label{initialarea}
\end{equation}
The scalars $r$ and $d$ can arbitrarily be chosen. Geometrical properties depend  only on ratios of lengths, without changing the ratio itself  and thus without loss of generality $d=1$, so that the eq. \eqref{initialarea} takes the form
\begin{equation}
	  A_{0}=\int _0^1\int _0^1\sqrt{4 r^2 u^2-4 r^2 u+4 r^2 v^2-4 r^2 v+2 r^2+1} \hspace{0.2cm}du dv.
\end{equation}
Substituting $H_{0}$ from eq. \eqref{mc0} in eq. \eqref{vp}, we have
\begin{equation}
\begin{split}
	m(u,v) &= 16 r^3 u^3 v^3-24 r^3 u^3 v^2+8 r^3 u^3 v-24 r^3 u^2 v^3+36 r^3 u^2 v^2 \\ & \qquad -12 r^3 u^2 v+8 r^3 u v^3-12 r^3 u v^2+4 r^3 u v. \label{mc1}
\end{split}
\end{equation}
Using \eqref{mc1}, variationally improved surface eq. \eqref{VS1} takes following form
\begin{equation} \label{VS4L}
\begin{split}
\mathbf x_{1}(u,v,t) &=	(16 r^3 t u^3 v^3-24 r^3 t u^3 v^2+8 r^3 t u^3 v-24 r^3 t u^2 v^3+36 r^3 t u^2 v^2-12 r^3 t u^2 v \\ & \qquad +  8 r^3 t u v^3-12 r^3 t u v^2+4 r^3 t u v-2 r u v+r u+r v,v,u ).
\end{split}
\end{equation}
Fundamental magnitudes for this variationally improved surface \eqref{VS4L} are as follows:
\begin{align}
\begin{split}
	E_{1}(u,v,t) &= (-t (8 r^2 (1-u) u (1-v) v (r (1-v)-r v)-4 r (1-u) (1-v) v (r (1-u)-r u) (r (1-v)-r v) \\ & \qquad +4 r u (1-v) v (r (1-u)-r u) (r (1-v)-r v))+r (1-v)-r v )^2+1,
\end{split}
\end{align}
\begin{align}
\begin{split}
	F_{1}(u,v,t) &= (-t (8 r^2 (1-u) u (1-v) v (r (1-u)-r u)-4 r (1-u) u (1-v) (r (1-u)-r u)  \\ & \qquad (r (1-v)-r v)+4 r (1-u) u v (r (1-u)-r u) (r (1-v)-r v))+r (1-u)\\ & \qquad -r u ) (-t ( 8 r^2 (1-u) u (1-v) v (r (1-v)-r v)-4 r (1-u) (1-v) v (r (1-u)-r u) \\ & \qquad  (r (1-v)-r v)+4 r u (1-v) v (r (1-u)-r u) (r (1-v)-r v))+r (1-v)-r v),
\end{split}
\end{align}
\begin{align}
\begin{split}
	G_{1}(u,v,t) &= (-t (8 r^2 (1-u) u (1-v) v (r (1-u)-r u)-4 r (1-u) u (1-v) (r (1-u)-r u) \\ & \qquad (r (1-v)-r v) +4 r (1-u) u v (r (1-u)-r u) (r (1-v)-r v))+r (1-u)-r u )^2+1,
\end{split}
\end{align}
\begin{align}
	e_{1}(u,v,t) = t (-96 r^3 u v^3+144 r^3 u v^2-48 r^3 u v+48 r^3 v^3-72 r^3 v^2+24 r^3 v),
\end{align}
\begin{align}
	f_{1}(u,v,t) = t (-144 r^3 u^2 v^2+144 r^3 u^2 v-24 r^3 u^2+144 r^3 u v^2-144 r^3 u v+24 r^3 u-24 r^3 v^2+24 r^3 v-4 r^3)+2 r,
\end{align}
and
\begin{align}
	g_{1}(u,v,t) = t (-96 r^3 u^3 v+48 r^3 u^3+144 r^3 u^2 v-72 r^3 u^2-48 r^3 u v+24 r^3 u).
\end{align}
Inserting these values of fundamental magnitudes in eq. \eqref{MC1} we find the expression for $H_{1}(u,v,t)$ of   surface  \eqref{VS4L} as
\begin{align}
\begin{split}
		H_{1}(u,v,t) &= [-4 r^3 (2 u-1) (2 v-1)]+[8 r^3 (2 u-1) (2 v-1) (r^2 (u^2 (6 v-5) (6 v-1) +u (-36 (v-1) v-5)+5 (v-1) \\ & \qquad v+1) -3 (u^2+v^2) +3 (u+v))]t+  [-32 r^7 (2 u-1) (2 v-1)  (6 u^4 (2 (v-1)  v (18 (v-1) v+5)+1)  -12 u^3 \\ & \qquad (2 (v-1) v (18 (v-1) v+5)+1)+u^2 (2 (v-1) v (138 (v-1) v  +37)+7) +u (-2 (v-1)   v (30 (v-1) \\ & \qquad v+7)-1)+(v-1) v (6 (v-1) v+1))]t^{2}	 +[u (2 u-1)  (v-1) v (2 v-1)  (12 u^4 (3 (v-1)  v (12 (v-1) \\ & \qquad  v+5)+2) -24 u^3 (3 (v-1) v (12 (v-1) v+5)+2)  +3 u^2 (12 (v-1) v  (17 (v-1) v+7)+11) \\ & \qquad -9 u (4 (v-1) v (5 (v-1) v+2)+1)+3 (v-1) v (8 (v-1) v+3)+1)]t^3.
\end{split}
\end{align}
After performing the integrations mentioned in eq. \eqref{rms1}, the  mean square curvature $\mu_{1}^{2}(t)$ for $\mathbf x_{1}(u,v)$ becomes
\begin{equation}
\begin{split}
	\mu_{1}^{2}(t) &= \left(\frac{2048 r^{18}}{2277275}\right)t^6+\left(\frac{190464 r^{16}}{25050025} \right)t^5+\left(\frac{512 r^{12} \left(153 r^2+77\right)}{444675} \right)t^4+\left( \frac{256 r^{10} \left(7 r^2+3\right)}{3675}\right)t^3 \\ & \qquad +\left(\frac{32 r^6 \left(29 r^4+98 r^2+119\right)}{1225} \right)t^2+\left(-\frac{64}{75} r^6 \left(3 r^2+5\right) \right)t+\left(\frac{16 r^6}{9} \right),
\end{split}	
\end{equation}
which may be minimized for $t$ for every fixed value of $r$. Fig. \ref{interpolation1} represents this minimizing value of $t_{min}$ as the numerical function of $r$.

\begin{figure}[htb!]
\begin{center}
\includegraphics[width=60mm]{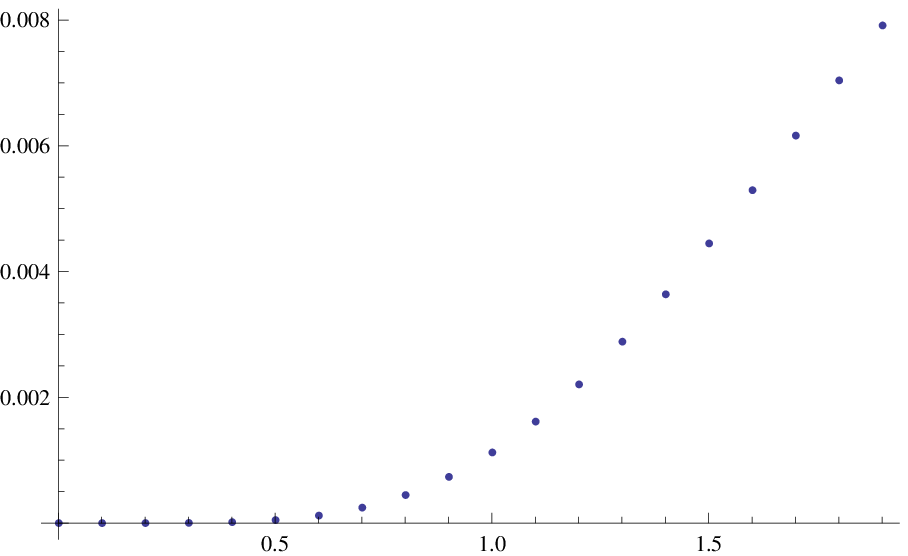}
\end{center}
\caption{The dimensionless decrease in area $A(r)$ as a function of $r$ enclosed by a skew quadrilateral ${\text {ruled}}_1$  bounded by  four arbitrary straight lines connecting four corners $ \mathbf x(0,0),\mathbf  x(0,1), \mathbf  x(1,0)$ and $\mathbf  x(1,1)$.}\label{listplot2}
\end{figure}

\begin{figure}[htb!]
\begin{center}
\includegraphics[width=60mm]{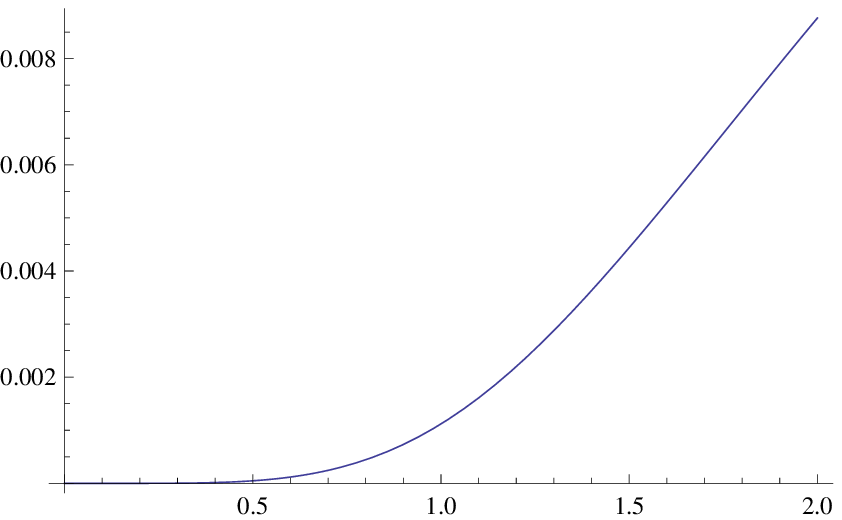}
\end{center}
\caption{The dimensionless decrease in area $A(r)$ as a function of  $r$ enclosed by a skew quadrilateral ${\text {ruled}}_1$  bounded by  four arbitrary straight lines
connecting four corners  $\mathbf  x(0,0), \mathbf x(0,1), \mathbf  x(1,0)$ and $\mathbf  x(1,1)$.} \label{interpolation2}
\end{figure}

We find the variationally improved surface $\mathbf x_{1}\left(u,v \right)$ eq. \eqref{VS1}  and its area as given by eq.  \eqref{area1}  for each  $t_{min}$  for the corresponding  $r$.  For a selection of $r$ values for $0 \leq r \leq 2$ with step size $0.001$, the  dimension less decrease in area  of surface $\mathbf x\left(u,v \right)$ of eq. \eqref{BIE} can be seen in the Fig. \ref{listplot2} and interpolating curve of the same is provided  in Fig.  \ref{interpolation2}.
\section{Conclusions} \label{conclusion}
We have discussed a technique to reduce the area of a surface $\mathbf x(u,v)$  eq. \eqref{CPR} obtaining variationally improved surface $\mathbf x_{1}(u,v)$ of eq. \eqref{VS1}.  This algorithm is first applied to a non- minimal surface spanning a boundary for which minimal surface is known namely hemiellipsoid  eq. \eqref{hemiellipsoid}. The dimensionless decrease $p$  in the area of the  hemiellipsoid  eq. \eqref{hemiellipsoid} for different values of $b$ and $c$ is $0\leq p \leq 15$ (see Fig. \ref {hemiellipsoid3}) depending upon  how much it is far from the minimal surface,  namely the elliptic disk.  This shows our algorithm eq. \eqref{VS1}   can significantly reduce area of surface that is far from being minimal. After noting this effectiveness, we applied this technique to reduce the area of a surface  of   $\mathbf x(u,v)$ eq. \eqref{BIE} bilinearly spanned by four non-planar boundary lines, a special case of Coons patch eq. \eqref{CPR},    along with  the configuration eq. \eqref{ruled1},  for a selection of $r$ values for $0 \leq r \leq 2$ with step size $0.001$. This gave us a much lesser ( in the range  $0$ to $0.80$)  dimensionless decrease in less area of surface $\mathbf x\left(u,v \right)$ of eq. \eqref{BIE}, as seen in the Fig. \ref{listplot2} or Fig.  \ref{interpolation2}. This suggests that ${\text {ruled}}_1$   is already a near minimal surface.
\bibliographystyle{unsrt}
\bibliography{XBibP1V2}

\end{document}